\documentclass[12pt]{article} 
\usepackage{amsmath,amsthm, amssymb} 
\usepackage{amssymb,latexsym}

\newtheorem*{theorem}{Theorem}

\begin{document}

\baselineskip=17pt 

\title{The ternary Goldbach problem with primes from arithmetic progressions}

\author{D. I. Tolev \footnote{Supported by Sofia University Grant 221/2008}}

\date{}
\maketitle

In 1937 Vinogradov \cite{Vin} found an asymptotic formula for the number of
solutions of the equation
\begin{equation} \label{l1}
p_1 + p_2 + p_3 = N
\end{equation}
in prime numbers $p_1, p_2, p_3$. Suppose that $k_i, l_i$; $i=1,2,3$ are integers satisfying $(k_i, l_i)=1$.
Using Vinogradov's method Zulauf \cite{Zul} established an asymptotic formula for the number of solutions
of \eqref{l1} in primes $p_i \equiv l_i \pmod{k_i}$, $i=1,2,3$.
This formula is valid not only for fixed $k_i$, but also for $k_i \le \mathcal L^D$, where 
$\mathcal L = \log N$ and $D>0$ is a constant. More precisely, in this case we have
\begin{equation} \label{l2}
\mathcal R_{{\bf k}, {\bf l}}(N) := \sum_{ \substack{ p_1 + p_2 + p_3 = N \\
    p_i \equiv l_i \pmod{k_i} \\ i=1,2,3  } } \log p_1 \log p_2 \log p_3 \; = \;
     \mathcal M_{{\bf k}, {\bf l}}(N) +  O \left(N^2 \mathcal L^{-A} \right) ,
\end{equation} 
where $A>0$ is arbitrarily large constant. Here the main term is given by
\[
\mathcal M_{{\bf k}, {\bf l}}(N) = \frac{N^2 \, \mathfrak S_{{\bf k}, {\bf l}} (N)}{2 \varphi(k_1) \varphi(k_2)\varphi(k_3)} ,
\]
where $\mathfrak S_{{\bf k}, {\bf l}} (N)$ is the singular series (see \cite{Hal-3} for the explicit formula)
and $\varphi(k)$ stands for the Euler function.

\bigskip

One may expect that \eqref{l2} is valid even if one or more of the moduli 
$k_i$ grow faster than a power of $\mathcal L$ but this has not been proved so far.
One may try instead to prove that \eqref{l2} is true 
on average with respect to $k_i$ from a certain set of positive integers. More precisely, define
\[
\Delta_{{\bf k}, {\bf l}}(N) = \mathcal R_{{\bf k}, {\bf l}}(N) - \mathcal M_{{\bf k}, {\bf l}}(N) 
\]
and consider the sum
\[
\mathcal E = \mathcal E(H_1, H_2, H_3) =
\sum_{ \substack{ k_i \le H_i \\ i=1,2,3}} 
\max_{ \substack{ (l_i, k_i)= 1  \\ i=1,2,3 }  }
 \left| \Delta_{{\bf k}, {\bf l}}(N)  \right|
\]
Bearing in mind Bombieri-Vinogradov's theorem, one may expect 
that for any constant $A>0$ there exists $B=B(A)>0$ such that 
\begin{equation} \label{l3}
 \mathcal E(H_1, H_2, H_3) \ll N^2 \mathcal L^{-A} \qquad \text{for} \qquad
 H_1, H_2, H_3 \le \sqrt{N} \mathcal L^{-B} .
\end{equation}

\bigskip

During the last years several results that approximate this conjecture were found,
and some of them were applied for studying the equation \eqref{l1}, or similar ternary problems, with prime variables 
having certain additional properties. The author~\cite{Tol-1} proved a theorem of this type, but with only one prime lying in a progression. This result was later improved by K.Halupczok~\cite{Hal-2}. A theorem with two primes from independent progressions (and with fixed $l_i$) was established by Peneva and the author \cite{Pen-Tol}. A stronger result was recently found by K.Halupczok~\cite{Hal-3}. 
Theorem~1 of \cite{Hal-3} states (essentially) that for any $A>0$ there exists $B=B(A) > 0 $ such that 
\begin{equation} \label{l4}
\sum_{k_1 \le \sqrt{N} \mathcal L^{-B}} \max_{(l_1, k_1)=1} \sum_{k_2 \le \sqrt{N} \mathcal L^{-B}} \max_{(l_2, k_2)=1}
\left| \Delta_{\{k_1, k_2, 1\}, \{l_1, l_2, 1\}}(N)  \right| 
\ll 
  N^2 \mathcal L^{-A} .
\end{equation}

\bigskip

It is not known at present if the estimate \eqref{l3} is true for $H_1 = H_2 = H_3 = N^{\delta}$, where $\delta >0$ is a constant. However using author's method from \cite{Tol-2} one can prove that 
if $l_1, l_2, l_3$ are fixed integers and $\lambda_i(k)$, $i=1,2,3$, are real numbers satisfying
$|\lambda_i(k)| \le 1$ then for any $A>0$ there exist $B=B(A)>0$
such that
\[
\sum_{\substack{k_1, k_2 \le \sqrt{N} \mathcal L^{-B} \, , \; k_3 \le N^{1/3} \mathcal L^{-B} \\ (l_i, k_i)=1 
 \, , \; i=1,2,3  }} 
\lambda_1(k_1) \lambda_2(k_2) \lambda_3(k_3)
\Delta_{{\bf k}, {\bf l}}(N) \ll N^2 \mathcal L^{-A} .
\]
(The upper bound for $k_3$ in the last formula can be increased to $N^{4/9} \mathcal L^{-B}$
in the case when 
$\lambda_3(k)$ is a well-factorable function of level $N^{4/9} \mathcal L^{-B}$. 
This is a consequence of a theorem of Mikawa \cite{Mik}.)

\bigskip

In this paper we present a new result that improves the theorems mentioned above. 
We have the following:
\begin{theorem} Suppose that $l_3$ is a fixed positive integer and 
$\lambda(k)$ are real numbers satisfying $|\lambda (k)| \le 1$.
For any constant $A>0$ there exists $B=B(A)>0$ such that 
\begin{equation} \label{l5}
\mathcal E^* := 
\sum_{\substack{k_1 \le \sqrt{N} \mathcal L^{-B} \\ k_2 \le \sqrt{N} \mathcal L^{-B}}}
\max_{\substack{ (l_1, k_1)=1 \\ (l_2, k_2)=1 }}
\left| \sum_{\substack{ k_3 \le N^{1/3} \mathcal L^{-B} \\ (k_3, l_3)=1 }} \lambda (k_3) \;
\Delta_{{\bf k}, {\bf l}}(N)  \right| \ll   N^2 \mathcal L^{-A} .
\end{equation}
\end{theorem} 

\bigskip

In particular, we find a stronger version of \eqref{l4} (with the maximum over
$l_1$ inside the sum over $k_2$; the method of \cite{Hal-3} is not applicable for proving this). 
We also mention that
if we apply Mikawa's theorem from \cite{Mik} then
the upper bound for $k_3$ in formula \eqref{l5} can be increased to $N^{4/9} \mathcal L^{-B}$ 
in the case when $\lambda(k)$ is a well-factorable function of level $N^{4/9} \mathcal L^{-B}$. 

\bigskip

We use the common number-theoretic notations. 
By greek letters we denote real numbers and 
by small latin letters --- integers. However, the letter $p$, with or without subscripts, is reserved for primes.
$N$ is a sufficiently large odd integer and $\mathcal L = \log N$.
As usual $\tau(k)$ is the number of positive divisors of $k$. By $(k, l)$ we denote the greatest common divisor of $k$ and $l$. Instead of $m \equiv n \pmod{k}$ we write $m \equiv n \, (k)$
and $k \sim K$ is abbreviation of $K < k \le 2K$. We also denote
$e(\alpha) = \exp(2 \pi i \alpha)$ and $||\alpha|| = \min_{n \in \mathbb Z} |\alpha - n|$.

\bigskip

{\bf Proof of the Theorem:} We may assume that
\begin{equation} \label{l5.5}
   A \ge 10^4 , \qquad B = 10^4 A.
\end{equation}
We apply the circle method with
$Q = \mathcal L^{20A}$, $\tau = N Q^{-1}$ and with the sets of major arcs $\mathfrak M$ and minor arcs $\mathfrak m$ 
specified by
\[
\mathfrak M = \bigcup_{q \le Q} \bigcup_{\substack{a=0 \\ (a, q)=1}}^{q-1} 
\left(\frac{a}{q} - \frac{1}{q\tau} , \frac{a}{q} + \frac{1}{q\tau} \right) , \qquad
\mathfrak  m = \left( -\frac{1}{\tau} , 1 - \frac{1}{\tau} \right) \setminus \mathfrak  M .
\]
It is clear that
\[
\mathcal R_{ {\bf k}, {\bf l} }(N) = \int_0^1 S_{k_1, l_1} (\alpha) S_{k_2, l_2} (\alpha) S_{k_3, l_3} (\alpha) 
   e(-N \alpha)d \alpha 
   = \mathcal R^{(\mathfrak  M)}_{ {\bf k}, {\bf l} } + \mathcal R^{(\mathfrak m)}_{ {\bf k}, {\bf l} }   ,
\]
where $\mathcal R^{(\mathfrak  M)}_{ {\bf k}, {\bf l} }$ and $\mathcal  R^{(\mathfrak m)}_{ {\bf k}, {\bf l} } $ are respectively the contributions coming from the major arcs and the minor arks and where
\begin{equation} \label{l6}
S_{k, l} (\alpha) = \sum_{ \substack{p \le N \\ p \equiv l \, (k) }} (\log p) \, e(\alpha p) .
\end{equation}

\bigskip

We have
\begin{equation} \label{l7}
\mathcal E^* \ll \mathcal E_1 + \mathcal E_2 .
\end{equation}
where
\begin{align}
\mathcal E_1 
 & =
  \sum_{\substack{k_1, k_2 \le \sqrt{N} \mathcal L^{-B} \\ i=1,2 }}
\max_{\substack{ (l_i, k_i)=1 \\
    i=1,2 }}   \left| 
    \sum_{\substack{ k_3 \le N^{1/3} \mathcal L^{-B} \\ (k_3, l_3)=1 }} \lambda (k_3)  
  \left(
\mathcal R^{(\mathfrak M)}_{ {\bf k}, {\bf l} } - \mathcal M_{{\bf k}, {\bf l}}(N) \right)       
   \right| ,
  \notag \\
 \mathcal E_2 
   & =
  \sum_{\substack{k_1, k_2 \le \sqrt{N} \mathcal L^{-B} \\ i=1,2 }}
\max_{\substack{ (l_i, k_i)=1 \\
    i=1,2 }}  
    \left| 
         \sum_{\substack{ k_3 \le N^{1/3} \mathcal L^{-B} \\ (k_3, l_3)=1 }} \lambda (k_3)  
  \mathcal R^{(\mathfrak m)}_{ {\bf k}, {\bf l} } 
      \right| = 
      \sum_{\substack{k_1, k_2 \le \sqrt{N} \mathcal L^{-B} \\ i=1,2 }}
\max_{\substack{ (l_i, k_i)=1 \\
    i=1,2 }}  
    \left| \mathcal U  \right| ,
    \notag
\end{align}
say. 

\bigskip

Working as in sections 4 and 5 of \cite{Tol-2} (see also Theorem 3 of \cite{Hal-3}) we find
\begin{equation} \label{l8}
 \mathcal E_1 \ll  \sum_{\substack{k_i \le \sqrt{N} \mathcal L^{-B} \\ i=1,2,3 }} 
\max_{\substack{ (l_i, k_i)=1 \\
    i=1,2,3 }}   \left| 
  \mathcal R^{(\mathfrak M)}_{{\bf k}, {\bf l}} - \mathcal M_{{\bf k}, {\bf l}}(N)   \right|  \ll N^2 \mathcal L^{-A} .
\end{equation}

\bigskip

It remains to estimate $\mathcal E_2$. Obviously
\begin{equation} \label{l9}
\mathcal E_2 \ll \mathcal L^2 \max_{K_1, K_2 \le \sqrt{N} \mathcal L^{-B}}
\mathcal V (K_1, K_2) ,
\end{equation}
where
\begin{equation} \label{l9.5}
\mathcal V = \mathcal V (K_1, K_2) = \sum_{ \substack{ k_1 \sim K_1 \\ k_2 \sim K_2 } } \max_{ \substack{ (l_i, k_i)=1 \\
    i=1,2 } }  
\left|  \mathcal U \right| .
\end{equation}
Using the definitions of $\mathcal R^{(\mathfrak m)}_{ {\bf k}, {\bf l} }$ and $\mathcal U$ we get
\[
\mathcal U = \int_{\mathfrak m} S_{k_1, l_1} (\alpha) S_{k_2, l_2} (\alpha) \mathcal K (\alpha)
e(- N \alpha) d \alpha ,
\]
where
\[
\mathcal K (\alpha) = 
\sum_{\substack{ k_3 \le N^{1/3} \mathcal L^{-B} \\ (k_3, l_3)=1 }} \lambda (k_3) S_{k_3, l_3} (\alpha) .
\]

\bigskip

This sum behaves, in some sense, like $S(\alpha)= \sum_{p \le N} (\log p) e (\alpha p)$.
More precisely, the following estimates hold:
\begin{equation} \label{l10}
\max_{\alpha \in \mathfrak m} \left| \mathcal K (\alpha) \right| \ll N \mathcal L^{350 - 4A} ,
   \qquad \qquad \int_0^1 \left| \mathcal K (\alpha) \right|^2 d \alpha \ll N \mathcal L^{20} .
\end{equation}
The proof of the first one can be obtained 
following the proof of Lemma~12 of \cite{Tol-2}, whiles the verification of the second is simple (see \cite{Tol-2}, page 88).

\bigskip

Using \eqref{l6} we find
\[
\left| \mathcal U \right| =
\left| \sum_{\substack{ p \le N \\ p \equiv l_2 \, (k_2) }} (\log p) 
I_{k_1, l_1}(p) \right| \ll
\mathcal L \sum_{\substack{ r \le N \\ r \equiv l_2 \, (k_2) }} 
\left| I_{k_1, l_1}(r) \right| ,
\]
where we have denoted
\begin{equation} \label{l12}
I_{k, l}(r) = \int_{\mathfrak m} S_{k, l}(\alpha)  \mathcal K (\alpha) 
e((r-N) \alpha ) d \alpha .
\end{equation}
Applying Cauchy's inequality we get
\begin{align} 
\mathcal V  
  & 
  \ll  \mathcal L (K_1 K_2)^{1/2} \left(
  \sum_{ \substack{ k_1 \sim K_1 \\ k_2 \sim K_2 } } \max_{ l_1, l_2 } 
 \Big( \sum_{\substack{ r \le N \\ r \equiv l_2 \, (k_2) }} 
\big| I_{k_1, l_1}(r) \big| 
  \Big)^2 
  \right)^{1/2}
    \notag \\
  &
   \ll  \mathcal L 
  (N K_1)^{1/2}
  \left(   \sum_{ \substack{ k_1 \sim K_1 \\ k_2 \sim K_2 } } \max_{ l_1, l_2 } 
  \mathcal F \right)^{1/2} ,
  \label{l13}
\end{align}
where
\[
\mathcal F  = \sum_{\substack{ r \le N \\ r \equiv l_2 \, (k_2) }} 
\big| I_{k_1, l_1}(r) \big|^2 .
\]

\bigskip

Consider $\mathcal F$. First we use \eqref{l12} and expand the square and then we insert the summation over $r$ inside 
the double integral:
\begin{align}
\mathcal F
  & =
\sum_{\substack{ r \le N \\ r \equiv l_2 (k_2) }} 
\int_{\mathfrak m} \int_{\mathfrak m} S_{k_1, l_1}(\alpha) \mathcal K (\alpha) S_{k_1, l_1}(-\beta) \mathcal K (-\beta)
e((r-N)(\alpha - \beta)) d \alpha d \beta
 \notag \\
  & =
  \int_{\mathfrak m} \int_{\mathfrak m} S_{k_1, l_1}(\alpha) \mathcal K (\alpha) S_{k_1, l_1}(-\beta) \mathcal K (-\beta)
  \sum_{\substack{ r \le N \\ r \equiv l_2 (k_2) }}
  e((r-N)(\alpha - \beta)) d \alpha d \beta .
\notag 
\end{align}
Now we estimate the sum over $r$ using the well-known bound for the linear exponential sum
and then apply the inequality $uv \le u^2+v^2$ to get
\begin{align}
\mathcal F
  & \ll
  \int_{\mathfrak m} \int_{\mathfrak m} \left| S_{k_1, l_1}(\alpha) \mathcal K (\alpha) S_{k_1, l_1}(\beta) \mathcal K (\beta) \right|
  \min \left(
     \frac{N}{K_2} , \frac{1}{||(\alpha - \beta) k_2  ||}
        \right)
   d \alpha d \beta
   \notag  \\
   & \ll
   \int_{\mathfrak m} \int_{\mathfrak m} \left| S_{k_1, l_1}(\alpha)  \mathcal K (\beta) \right|^2
  \min \left(
     \frac{N}{K_2} , \frac{1}{||(\alpha - \beta) k_2  ||}
        \right)
   d \alpha d \beta  .
    \notag  
\end{align}

\bigskip

Next we extend the integration over $\alpha$ to the unit interval and use that
\[
|S_{k_1, l_1}(\alpha)|^2 = \sum_{\substack { |n| \le N \\ n \equiv 0 \, (k_1) }} w(n) \, e(\alpha n) ,
\]
where
\begin{equation} \label{l14}
w(n) = w(n, k_1, l_1) = \sum_{\substack{ p_1, p_2 \le N \\ p_1 - p_2 = n \\ p_1 \equiv l_1 \, (k_1) }}
\log p_1 \log p_2  \ll \mathcal L^2 N K_1^{-1} .
\end{equation}
We obtain
\begin{align}
\mathcal F
   & \ll
   \int_{\mathfrak  m } |\mathcal K (\beta)|^2
   \int_0^1 
   \sum_{\substack { |n| \le N \\ n \equiv 0 \, (k_1) }} w(n) 
   e(n \alpha)
     \min \left(
     \frac{N}{K_2} , \frac{1}{||(\alpha - \beta) k_2  ||}
        \right)   d \alpha  d \beta
    \notag \\
    & \ll
    \sum_{\substack { |n| \le N \\ n \equiv 0 \, (k_1) }} |w(n)|
    \left|
    \int_{\mathfrak  m } |\mathcal K (\beta)|^2
   \int_0^1 
    e(n \alpha)
     \min \left(
     \frac{N}{K_2} , \frac{1}{||(\alpha - \beta) k_2  ||}
        \right)   d \alpha  d \beta
        \right| .
    \notag        
\end{align}

\bigskip

We use \eqref{l14} and change the variable in the inner integral to get
\begin{align}
\mathcal F 
  & \ll 
  \frac{\mathcal L^2 N}{K_1} 
\sum_{\substack { |n| \le N \\ n \equiv 0 \, (k_1) }}
  \left|
  \int_{ \mathfrak m } |\mathcal K (\beta)|^2
  \int_{-\beta}^{1-\beta} 
    e(n (\beta + \gamma))
     \min \left(
     \frac{N}{K_2} , \frac{1}{||\gamma k_2  ||} \right)
     d \gamma  d \beta
     \right|
     \notag \\
  & =
   \frac{\mathcal L^2 N}{K_1} 
\sum_{\substack { |n| \le N \\ n \equiv 0 \, (k_1) }}
  \left|
  \int_{\mathfrak  m } |\mathcal K (\beta)|^2 e(n \beta) d \beta \right|
   \; \left| J(n, k_2)
       \right| ,
     \notag
\end{align}
where 
\[
J(n, k) =
\int_0^1 
    e(n \gamma)
     \min \left(
     \frac{N}{K_2} , \frac{1}{||\gamma k  ||} \right)
     d \gamma  .
\]     
To study this integral we apply the well-known decomposition
\[
 \min \left( H, \frac{1} {|| \alpha||} \right) = \sum_{h} c(h) e(h \alpha) ,
 \qquad
 c(h) \ll
 \begin{cases}
   \log H       &  \text{ for all } h , \\
   H^2 h^{-2}   &  \text{ for } h \not= 0 ,
  \end{cases}
\]
(with $H = N/K_2$) and we find that
\begin{align}
J(n, k) 
  & = 
  \int_0^1 e(n \gamma) \sum_{h} c(h) e(h k \gamma) d \gamma =
 \sum_h c(h) \int_0^1 e ( (n + hk) \gamma ) d \gamma
 \notag \\
 &  = 
 \begin{cases} 
 c(-n/k)  & \text{ if } k \mid n , \\
 0        & \text{ if } k \nmid n .
 \end{cases}
 \notag 
\end{align}

\bigskip

Therefore
\begin{equation} \label{l15}
\mathcal F \ll
\mathcal L^3 N K_1^{-1}
\sum_{\substack { |n| \le N \\ n \equiv 0 \, (k_1) \\ n \equiv 0 \, (k_2) }}
  \left|
  \int_{ \mathfrak m } |\mathcal K (\beta)|^2 e(n \beta) d \beta \right| .
\end{equation}
We note that the last expression already does not depend on $l_1$ and $l_2$.
Using \eqref{l13} and \eqref{l15} we get
\begin{equation} \label{l16}
\mathcal V \ll
  \mathcal L^3 N
  \left( 
  \sum_{ \substack{ k_1 \sim K_1 \\ k_2 \sim K_2 } } 
\sum_{\substack { |n| \le N \\ n \equiv 0 \, (k_1) \\ n \equiv 0 \, (k_2) }}
  \left|
  \int_{\mathfrak  m } |\mathcal K (\beta)|^2 e(n \beta) d \beta \right| \right)^{1/2} =
  \mathcal L^3 N \, \mathcal W^{1/2} ,
\end{equation}
say. 

\bigskip

Now we write
\begin{equation} \label{l17}
\mathcal W = \mathcal W' + \mathcal W^* ,
\end{equation}
where $\mathcal W'$ is the contribution of the terms with $n=0$ and $\mathcal W^*$ comes from the other terms.
From \eqref{l9} and \eqref{l10} we get
\begin{equation} \label{l18}
\mathcal W' \ll K_1 K_2 \int_0^1 |\mathcal K (\beta)|^2  d \beta  \ll 
  \mathcal L^{20} K_1 K_2 N \ll \mathcal L^{20 - 2 B} N^2 .
\end{equation}

\bigskip

Consider $\mathcal W^*$. We have
\[
\mathcal W^*  =
   \sum_{1 \le |n| \le N} \sum_{\substack{ k_i \sim K_i  \\ k_i \mid n , \, i=1,2  }}
  \left|
  \int_{\mathfrak  m } |\mathcal K (\beta)|^2 e(n \beta) d \beta \right| 
  \ll 
  \sum_{1 \le n \le N} \tau^2(n) 
  \left|
  \int_{\mathfrak  m } |\mathcal K (\beta)|^2 e(n \beta) d \beta \right| .
\]
Applying the inequalities of Cauchy and Bessel and using the well-known elementary estimate
$\sum_{n \le N} \tau^4(n) \ll N \mathcal L^{15}$ we find that
\begin{align}
\mathcal W^*
  & \ll
  \left( \sum_{ n \le N} \tau^4(n) \right)^{1/2}
      \left( \sum_{ n \le N} \left|
  \int_{\mathfrak  m } |\mathcal K (\beta)|^2 e(n \beta) d \beta \right|^2
  \right)^{1/2} 
  \notag \\
  & \ll
  \mathcal L^8 N^{1/2} 
   \left(   \int_{\mathfrak  m } |\mathcal K (\beta)|^4  d \beta  \right)^{1/2} 
   \notag \\
   & \ll 
   \mathcal L^8 N^{1/2}   \sup_{\beta \in \mathfrak m} |\mathcal K (\beta)| 
      \left(
   \int_0^1 |\mathcal K (\beta)|^2  d \beta \right)^{1/2} .
   \notag
     \end{align}
We estimate the last expression using  \eqref{l10} to obtain
\begin{equation} \label{l19}
\mathcal W^* \ll N^2 \mathcal L^{400 - 4A} .
\end{equation}

From \eqref{l5.5}, \eqref{l9} and \eqref{l16} -- \eqref{l19} we get
\[
\mathcal E_2 \ll N^2 \mathcal L^{-A} 
\]
and the proof of the theorem follows from this estimate, \eqref{l7} and \eqref{l8}.

\bigskip
\bigskip

\vbox{
\hbox{Faculty of Mathematics and Informatics}
\hbox{Sofia University ``St. Kl. Ohridsky''}
\hbox{5 J.Bourchier, 1164 Sofia, Bulgaria}
\hbox{ }
\hbox{Email: dtolev@fmi.uni-sofia.bg}}


\begin{thebibliography}{99}


\bibitem{Hal-2}
K. Halupczok, {\it On the number of representations in the ternary Goldbach problem with one prime number in a given residue class}, J. Number Theory 117 (2006), no. 2, 292--300. 

\bibitem{Hal-3}
K. Halupczok, {\it On the ternary Goldbach problem with primes in independent arithmetic progressions}, Acta Math. Hungar., 120 (4) (2008), 315-349. 


\bibitem{Mik} H. Mikawa, {\it On exponential sums over primes in arithmetic progressions}
Tsukuba J. Math. 24, 2, (2000), 351-360.

\bibitem{Pen-Tol} T. P. Peneva and D. I. Tolev, {\it An additive problem with primes and almost-primes},
Acta Arith. 83 (1998), 155–169.


\bibitem{Tol-1} D. I. Tolev, {\it On the number of representations of an odd integer as a sum of three
primes, one of which belongs to an arithmetic progression}, Proc. Steklov. Inst. Math.
218 (1997).


\bibitem{Tol-2} D. I. Tolev, {\it Arithmetic progressions of prime-almost-prime
twins}, Acta Arith., 88, (1999), 67-98.


\bibitem{Vin} I. M. Vinogradov, {\it Representation of an odd number as the sum of three primes},
Dokl. Akad. Nauk. SSSR, 15, (1937), 291-294, (in Russian).


\bibitem{Zul} A Zulauf, 
{\it \"Uber die Darstellung nat\"urlicher Zahlen als Summen von Primzahlen aus gegebenen Restklassen und Quadraten mit gegebenen Koeffizienten, I: Resultate f\"ur gen\"ugend grosse Zahlen}
J. Reine Angew. Math. 192, (1954), 210-229;
{\it II:  Die singul\"are Reihe}, J. Reine Angew. Math. 193 (1954), 39-53. 




\end{thebibliography}
\end{document}